\documentclass[12pt, reqno]{amsart}

\newcommand{\tsum}{\textstyle\sum}

\newcommand{\boxt}{\Tox\kern -6.3pt\raise .55pt
    \hbox{$\scriptstyle{\times}$}}%boxtimes

\newcommand{\bg}{{\mathbf g}} 
\newcommand{\bh}{{\mathbf h}}

\newcommand{\cdott}{\cdot,\cdot}

\newcommand{\ltimes}{\vbox to 5.4pt{\leaders\vrule\vfil}\kern
  -5pt\times}

\newcommand{\del}{\partial}

\newcommand{\pdb}[2][{}]
          {\frac{\partial^{#1}\phantom{x}}{\partial#2}}

\newtheorem{thm}{Theorem}[section]
\newtheorem{lem}[thm]{Lemma}

\newtheorem{prop}[thm]{Proposition}
\newtheorem{cor}[thm]{Corollary}

\theoremstyle{definition}

\newtheorem{remark}[thm]{Remark}

\newtheorem{remarks}[thm]{Remarks}

\numberwithin{equation}{section}

\newcommand{\inv}{^{-1}} \newcommand{\invv}[1][]{\frac{1}{#1}}
 
\newcommand{\wt}{\widetilde}      

\newcommand{\half}[1][1]{\frac{#1}{2}}

          \newcommand{\g}{{\mathfrak g}}

\newcommand{\fsl}{\mathfrak{sl}}

\newcommand{\C}{{\mathbb C}}

\newcommand{\R}{{\mathbb R}}

\newcommand{\al}{\alpha}
\newcommand{\be}{\beta}

\newcommand{\la}{\lambda}

\newcommand{\ze}{\zeta}

                                     \newcommand{\A}{{\mathcal A}}
                                     \newcommand{\B}{{\mathcal B}}

\newcommand{\E}{{\mathcal E}}

                                         \newcommand{\K}{{\mathcal K}}

\newcommand{\cP}{{\mathcal P}} 
                                            \newcommand{\Q}{{\mathcal Q}} 
 
\newcommand{\cS}{{\mathcal S}} 
                                           
                                           \newcommand{\U}{{\mathcal U}}

\newcommand{\AND}{\qquad\mbox{and}\qquad}

\newcommand{\Pol}{\mathop{\mathrm{Pol}}\nolimits}

\newcommand{\gr}{\mathop{\mathrm{gr}}\nolimits}

\newcommand{\Ug}[1][]{{\mathcal U}_{#1}({\mathfrak g})}

\newcommand{\sal}{^{\al}}
\newcommand{\sbe}{^{\be}}

\newcommand{\rep}{representation}  
\newcommand{\reps}{representations}

\newcommand{\D}{{\mathfrak D}} \newcommand{\Dla}{\D^{\la}}

\newcommand{\RR}[1][]{{\mathcal R^{#1}}}
 
\newcommand{\CP}{\mathbb{CP}}
\newcommand{\RP}{\mathbb{RP}}

\newcommand{\smhalf}{{\scriptscriptstyle \frac{1}{2}}}

\newcommand{\DD}[1][\smhalf]{{\mathcal D}^{#1}}
\newcommand{\SLR}[1][n+1]{SL_{#1}(\R)}  
\newcommand{\SLC}[1][n+1]{SL_{#1}(\C)}
\newcommand{\slR}[1][n+1]{\fsl_{#1}(\R)}  
\newcommand{\slC}[1][n+1]{\fsl_{#1}(\C)}

\newcommand{\sm}{_{\infty}}
\newcommand{\PolTRP}[1][]{\Pol\sm^{#1}(T^*\,\RP^n)}
\newcommand{\DRP}[1][\smhalf]{\D^{#1}\sm(\RP^n)}
\newcommand{\TR}{T^*\R^n} 
\newcommand{\TRP}{T^*\RP^n} 
\newcommand{\TCP}{T^*\CP^n} \newcommand{\RTCP}{R(T^*\CP^n)}
\newcommand{\RPn}{\RP^{n}}

\newcommand{\opluss}[1][d]{\oplus_{#1=0}^{\infty}}
\newcommand{\Bla}{\B^{\la}}  
 
\newcommand{\sla}{_{\la}}     
\newcommand{\slat}[1][t]{_{\la;#1}}  

\newcommand{\starla}{\star_{\la}}

\newcommand{\Cla}[1][\la]{C^{#1}}
\newcommand{\etala}{\eta_{\la}}

\newcommand{\Csm}{\C_{\infty}}

\newcommand{\summo}[1][d]{\sum_{#1=1}^\infty\,}

\begin{document}   
 
\title 
{Non-Locality of Equivariant Star Products on $\mathbf{T}^*(\RPn)$}
\author{Ranee Brylinski}
\address{Department of Mathematics,
        Penn State University, University Park 16802}
\email{rkb@math.psu.edu}
\urladdr{www.math.psu.edu/rkb}
%\thanks{}
%\subjclass{??22E46, 17B35, 53D55, 17C20, 43A85} 
%\keywords{}
 
%\date{}

\begin{abstract}
Lecomte and Ovsienko constructed $SL_{n+1}(\R)$-equivariant 
quantization maps  $\Q_{\lambda}$ for   symbols   of differential   
operators on $\lambda$-densities on $\RP^n$. 

We derive some formulas for the  associated graded 
equivariant star products $\star_{\lambda}$ on  
the symbol algebra $\Pol(T^*\RP^n)$. These
give some measure of the failure of locality. 

Our main result expresses (for $n$ odd)
the coefficients $C_p(\cdot,\cdot)$ of 
$\star_{\lambda}$ when $\lambda=\half$  in terms of  some new 
$SL_{n+1}(\C)$-invariant algebraic bidifferential operators 
$Z_p(\cdot,\cdot)$  on $T^*\CP^n$
and the operators $(E+\frac{n}{2}\pm s)^{-1}$ where $E$ is the  
fiberwise Euler vector  field and $s\in\{1,2,\cdots,[\frac{p}{2}]\}$.
\end{abstract} 
  
\maketitle

\section{Introduction} 
\label{sec_1}      
Lecomte and Ovsienko (\cite{L-O}) constructed    $SL_{n+1}(\R)$-equivariant
quantization maps  $\Q\sla$ for symbols of  differential   
operators on $\lambda$-densities on $\RP^n$. 

We derive some formulas for the  associated graded 
equivariant star products 
${\phi\starla\psi}=\phi\psi+\tsum_{p=1}^\infty\Cla_p(\phi,\psi)t^p$ on  the
symbol algebra $\Pol\sm(\TRP)$. 
The star products  $\starla$ is  ``algebraic" in
that (Proposition \ref{prop:3}) it restricts to the   subalgebra  $\RR$ generated by
the momentum functions  $\mu^x$, $x\in\slR$. 

We compute some special values of 
$\phi\starla\psi$ in Proposition \ref{prop:4}.
We conclude in  Corollary \ref{cor:4}  that  
$\Cla_p(\cdott)$ fails to be bidifferential, except if $\la=\half$ and $p=1$.
The reason is that $\Cla_p(\cdott)$ involves operators of the form $(E+r)\inv$
where  $E$ is the   fiberwise Euler vector    field on $\TRP$ and $r$ is a positive
number.
  
In our main result (Theorem \ref{thm:5}), we write, for $n$ odd,
the coefficients $\Cla_p(\cdot,\cdot)$ when $\lambda=\half$ in terms of  
some new   $SL_{n+1}(\C)$-invariant algebraic bidifferential operators 
$Z_p(\cdot,\cdot)$  on $\CP^n$
and the operators $(E+\frac{n}{2}\pm s)^{-1}$ where
$s\in\{1,2,\cdots,[\frac{p}{2}]\}$. Our proofs in \S\ref{sec_4}-\S\ref{sec_5} are
applications of  the formulas  in \cite[\S4.5]{L-O} for $\Q\sla$.

The operator  $Z_p(\cdott)$ ($p\ge 2$) is quite subtle   as it has total
homogeneous degree $-p$. It is  not   the $p$th power of the Poisson tensor
(with respect to some coordinates)
because we can show that the total order of $Z_p(\cdott)$ is too large.
It would be very interesting to find a way to construct    
$Z_p$ using the method of Levasseur and Stafford (\cite{L-S}).

I thank Christian Duval and  Valentin Ovsienko for  several interesting
discussions.
   
\section{The Lecomte-Ovsienko quantization maps}  
\label{sec_2}    

In \cite{L-O},  Lecomte and Ovsienko  constructed, for each $\la\in\C$,  
an $\SLR$-equivariant (complex 
linear) quantization map $\Q\sla$ from $\A=\PolTRP$ to $\Bla=\DRP[\la]$.
Here $\A=\opluss\A^d$ is the graded Poisson algebra of smooth complex-valued
functions on $\TRP$ which are polynomial   along the cotangent fibers, and  
$\Bla=\cup_{d=0}^\infty\Bla_d$ is the filtered
algebra of smooth (linear) differential operators on   $\la$-densities on $\RPn$.
Then  $\Q\sla$ is a quantization map in the sense that  $\Q\sla$ is a vector space 
isomorphism and $\phi$ is the principal symbol of $\Q\sla(\phi)$ if $\phi\in\A^d$.

The natural action of $\SLR$   on $\RPn$   lifts canonically to  a Hamiltonian action
on $\TRP$ with moment map $\mu:\TRP\to\slR^*$.  The density line bundle on 
$\RPn$ is  homogeneous for $\SLR$. This geometry produces   natural (complex
linear) \reps\ of   $\SLR$  on  $\A$ and $\Bla$;  $\Q\sla$ is equivariant for these
\reps. 
  
The procedure of Lecomte and Ovsienko  was to construct  
(\cite[Thm.  4.1]{L-O})      an 
$\slR$-equivariant quantization map $\Q\sla$   from $\Pol\sm(\TR)$ to 
$\Dla\sm(\R^n)$, where $\R^n$ is the   big cell in
$\RPn$.  They show their map is  unique.
Then $\Q\sla$ restricts to a quantization map from $\A$ to $\Bla$ 
(\cite[Cor. 8.1]{L-O}).

We can represent points in $\RPn$ in homogeneous coordinates
$[u_0,\dots,u_n]$. Then $u_1,\dots,u_n$ are linear coordinates on the big cell 
$\R^n$ defined by $u_0=1$. These, together with the  conjugate momenta
$\xi_1,\dots,\xi_n$,  give Darboux coordinates on $\TR$.

For any  vector field $\eta$ on $\RPn$, let $\mu_\eta\in\A^1$ be its 
principal symbol and let $\etala$ be its Lie derivative  acting on
$\la$-densities so that  $\etala\in\Bla_1$.  
Then $\Q\sla(\mu_\eta)=\etala$; this follows   by \cite[\S4.3]{L-O}.

The quantization map $\Q\sla$ defines a  star product; see \cite[\S8.2]{L-O}. 
For  $\phi,\psi\in\A$,   we put
$\phi\starla\psi=\Q\slat\inv(\Q\slat(\phi)\Q\slat(\psi))$ where 
$\Q\slat$ is  the linear map $\A\to\Bla[t]$ such that
$\Q\slat(\phi)=t^d\Q\sla(\phi)$ if $\phi\in\A^d$.
Then $\starla$ makes $\A[t]$ into an associative algebra over  $\C[t]$.
This satisfies
\begin{equation}\label{eq:starla} 
\phi\starla\psi=\tsum_{p=0}^\infty\Cla_p(\phi,\psi)t^p
\end{equation}
where $\Cla_0(\phi,\psi)=\phi\psi$ and 
$\Cla_1(\phi,\psi)-\Cla_1(\psi,\phi)=\{\phi,\psi\}$.
Also $\Cla_p(\phi,\psi)\in\A^{j+k-p}$ if $\phi\in\A^j$ and $\psi\in\A^k$.
So $\starla$ is a graded star product  on $\A$.

We say that $\starla$ has \emph{parity} if
$\Cla_p(\phi,\psi)=(-1)^p\Cla_p(\psi,\phi)$; then
$\Cla_1(\phi,\psi)=\half\{\phi,\psi\}$.  

\begin{lem}\label{lem:2} 
$\starla$ has parity iff $\la=\half$.
\end{lem} 
\begin{proof} 
Let $\be:\Bla\to\B^{1-\la}$  be the canonical algebra 
anti-isomorphism  and let  $\al:\A\to\A$   be the Poisson
algebra anti-involution defined by
$\phi\sal=(-1)^d\phi$ if $\phi\in\A^d$. Then
$\Q\sla(\phi\sal)\sbe=\Q_{1-\la}(\phi)$ by \cite[Lem. 6.5]{L-O}. 
This implies $\Cla_p(\phi,\psi)=(-1)^p\Cla[1-\la]_p(\psi,\phi)$.
So we have parity if  $\la=\half$.
Otherwise parity is violated, already for $\Cla_1$. Indeed, if  
$\phi\in\A^0$ and $\mu\in\A^1$, then
$\phi\starla\mu=\phi\mu+\la\{\phi,\mu\}t$, and so
$\Cla_1(\phi,\mu)=\la\{\phi,\mu\}$ while
$\Cla_1(\mu,\phi)=-\Cla[1-\la]_1(\phi,\mu)=(\la-1)\{\phi,\mu\}$.
\end{proof}

\section{Algebraicity of $\starla$}  
\label{sec_3}    

Each $x\in\slR$ defines a vector field $\eta^x$ on $\TRP$. The principal symbols
$\mu^x=\mu_{\eta^x}$ are the momentum functions for $\SLR$.
The $\SLR$-equivariance of $\Q\sla$ is equivalent to $\slR$-equivariance, i.e.,
$\Q\sla(\{\mu^x,\phi\})= [\etala^x,\Q\sla(\phi)]$.
Then  $\Q\sla$ is $\slC$-equivariant, where 
we define $\mu^x$ and $\etala^x$ for $x\in\slC$ by
$\mu^{x+iy}=\mu^x+i\mu^y$ and so on.

The algebra $\RTCP$ of regular functions (in the sense of algebraic geometry) on
(the  quasi-projective complex algebraic variety)  $\TCP$ identifies, by restriction,
with a subalgebra $\RR$ of $\A$. Similarly the algebra of
$\Dla(\CP^n)$ of twisted algebraic (linear) differential operators for the 
formal  $\la$th power of the canonical bundle $\K$ identifies with a subalgebra
$\DD[\la]$ of  $\Bla$.

Then $\RR$ is generated by the momentum functions
$\mu^x$,  $\DD[\la]$ is generated by the operators $\etala^x$, and
$\gr\DD[\la]=\RR$. These statements follow, for instance, by 
\cite[Lem. 1.4 and Thm. 5.6]{Bo-Br}, since the proofs of the relevant results
there generalize immediately to the twisted case.
We get natural identifications 
$\RR=\cS/I$ and $\DD[\la]=\Ug/J$ where $I$ is graded 
Poisson  ideal in the symmetric algebra  $\cS=S(\slC)$, $J$ is a two-sided
ideal in  the enveloping algebra $\U=\U(\slC)$, and $\gr J=I$.

Notice  $\RR$ carries a natural \rep\ of $\SLC$, which then extends the
$\SLR$-symmetry it inherits from $\A$. 
\begin{prop}\label{prop:3} 
For every $\la$,  $\starla$ restricts to a graded
$G$-equivariant star product on the momentum algebra $\RR$.
\end{prop} 
\begin{proof} 
It suffices to check that $\Q\sla$ maps $\RR$ onto $\DD[\la]$
(which is stated for  $\la=0$ in \cite[\S1.5, Remark (c)]{L-O}). This follows 
easily in any number of ways. For instance,
the formula  for $\Q\sla$ in \cite[(4.15)]{L-O} implies
$\Q\sla(\xi_1^{a_1}\cdots\xi_n^{a_n})=
\frac{\del^{a_1}}{\del u_1^{a_1}}\cdots\frac{\del^{a_n}}{\del u_n^{a_n}}$. 
But $\{\xi_n^{d}\}_{d=0}^\infty$ and 
$\{\frac{\del^{d}}{\del u_n^{d}}\}_{d=0}^\infty$ are
complete sets of lowest weight vectors in $\RR$ and $\DD[\la]$.
\end{proof}  
  
\begin{remark}\label{rem:par_RR}  
The restriction of $\starla$ to $\RR$ has parity iff (i) $\la=\half$  or (ii)
$n=1$; see \cite[\S3]{A-B:starmin}. Notice that (ii) does not contradict the proof
of Lemma  \ref{lem:2}, as $\RR^0=\C$. 
\end{remark}

%\begin{remark}\label{rem:axioms}  
%We note that axioms (ii)-(iv) and (iii) in \cite{me:edqcot}  are satisfied by 
%$\Q\sla$. Indeed, (ii) and  (iii) follow by uniqueness, while (iv) follows by the  
%explicit construction in \cite{L-O}.  
%\end{remark}

\section{Some special values of  $\phi\starla\psi$}  
\label{sec_4}

$\Pol\sm(\TR)$ is the tensor product of two maximal Poisson commutative
subalgebras, namely the algebra 
$\Csm[u]=\Csm[u_1,\dots,u_n]$ of smooth functions on the big cell $\R^n$ and
and the polynomial algebra $\C[\xi]=\C[\xi_1,\dots,\xi_n]$.
Let $E$ be the fiberwise Euler vector field $\sum_{i=1}^n\xi_i\pdb{\xi_i}$.
Set  $D={\sum_{i=1}^n\frac{\partial^{2}}{\partial u_i\partial\xi_i}}$. 
\begin{prop}\label{prop:4} 
If $\phi\in\Csm[u]$  and  $\psi\in\C[\xi]$ then
$\phi\starla\psi=\bg\sla(\phi\psi)$ where  
\begin{equation}\label{eq:bg} 
\bg\sla=1+\summo[d]g\slat[d]D^dt^d\AND
g\slat[d]=\invv[d!]\prod_{j=0}^{d-1}\frac{-E-j-\la(n+1)}{2E+j+ n+1} 
\end{equation}
\end{prop}
\begin{proof}  
Let  $\Q_{norm}:\Pol\sm(\TR)\to\D\sm(\R^n)$ be the  normal ordering
quantization  map. The    construction of $\Q\sla$ in \cite{L-O} gives
$\Q\sla=\Q_{norm}\bh\sla$  where 
$\bh\sla=1+\summo[d]h\slat[d]D^dt^d$ and  $h\slat[d]$ are certain operators.
Here $\D\sm(\R^n)$  identifies with $\Dla\sm(\R^n)$ in the usual way.

In \cite[Th. 4.1]{D-L-O}  they give a very nice formula 
for the $h\slat[d]$ when  $\la=\half$.  Going back to  
\cite[(4.15)]{L-O}, we get a similar formula for all $\la$.  We find
%$\ovl[C]^k_{k-d}$ coresponds to $b_d(E)=\frac{E+\la(n+1)}{d(2E+d+n)}$
%acting on $\A^{k-d}$. Then the coefficient of $D^d$ in $\bh\sla$ is
%$b_d(E)b_{d-1}(E+1)\cdotsb_1(E+d-1)$. So
\begin{equation}\label{eq:bh} 
h\slat[d]=\invv[d!]\prod_{j=0}^{d-1}\frac{E+j+\la(n+1)}{2E+j+n+d} 
\end{equation}

Thus for $\phi,\psi\in\Pol\sm(\TR)$ we  have 
\begin{equation}\label{eq:star=short} 
\phi\starla\psi=\bg\sla(\bh\sla(\phi)\#\bh\sla(\psi)) 
\end{equation}
where $\#$  denotes the graded star product   defined by $\Q_{norm}$    and
$\bg\sla=\bh\sla\inv$.  We   find, directly from (\ref{eq:bh}) or using
\cite[(4.10)]{L-O}, that $\bg\sla$ is given by (\ref{eq:bg}).  

We know  
$\phi\#\psi=\sum_{p=0}^\infty N_p(\phi,\psi)t^p$ where   
$N_k(\phi,\psi)=\frac{1}{k!}\sum_{\al\in\{1,\dots,n\}^k}
\frac{\del^k\phi}{\del\xi_\al}\frac{\del^k\psi}{\del u_\al}$.
Now, for $\phi\in\Csm[u]$  and  $\psi\in\C[\xi]$,  (\ref{eq:star=short}) gives
$\phi\starla\psi=\bg\sla(\phi\psi)$.
\end{proof}     
 
\begin{cor}\label{cor:4} 
None of the operators $\Cla_p$  \textup{(}$p\ge 1$, $\la\in\C$\textup{)} 
is bidifferential on $\TR$, 
with one exception: $2C^{\smhalf}_1$ is the Poisson bracket.
\end{cor}
\begin{proof} 
We just showed that 
$\Cla_p(\phi,\psi)=g\slat[p]D^p(\phi\psi)$ if $\phi\in\Csm[u]$  and 
$\psi\in\C[\xi]$. This implies, if $\Cla_p$ is bidifferential, that $g\slat[p]$ is a
differential operator on $\TR$. Looking at our expression for $g\slat[p]$, we
deduce  $E+j+\la(n+1)=E+\half[j]+\half(n+1)$ for $j=0,\dots,p-1$. But   this
forces   $p=1$ and $\la=\half$. By parity,  $C^{\smhalf}_1=\half\{\cdott\}$.
\end{proof}
 
The corollary contradicts the  claim in \cite[\S8.2]{L-O}. They  no doubt meant 
that for each pair $j,k$, the restricted map $\Cla_p:\A^j\times\A^k\to\A^{j+k-p}$
is given by some bidifferential operator.

\section{Coefficients $\Cla_p$ for  $\la=\half$}  
\label{sec_5}    

In this section, we set $\la=\half$ and suppress the corresponding super(sub)scripts.
We put $E'=E+\half[n]$ where $E$ is the fiberwise  Euler vector field  on $\TRP$.
See \cite{A-B:gq} for an interpretation of the shift $\half[n]$.
Let $[m]$ be the greatest integer not exceeding $m$.  

We put $T_p=\prod_{i=1}^{[\frac{p}{2}]}(E'+i)$ and
$S_p=\prod_{i=1}^{[\frac{p}{2}]}(E'-i)$. 
These are both invertible on   $\A$  if  $n$ is odd.
Our main result is 
\begin{thm} \label{thm:5} 
Assume $n$ is odd and let $p\ge 1$.  Then $C_p$  has the  form   
\begin{equation}\label{eq:Cp=} 
C_p(\phi,\psi)=\invv[T_p]  Z_p\left(\invv[S_p]\phi,\invv[S_p]\psi\right),
\qquad  \phi,\psi\in\A
\end{equation}
where   $Z_p$ is an $\SLR$-invariant bidifferential operator  on $\TRP$.

$Z_p$ is uniquely  determined by \textup{(\ref{eq:Cp=})}, even if we just take
$\phi,\psi\in\RR$.
Thus  $\star$ is uniquely determined by its restriction to $\RR$, once we know  
that  $(\phi,\psi)\mapsto T_pC_p(S_p\phi,S_p\psi)$ is bidifferential.

Finally, $Z_p$, like $E'$,  extends  uniquely to an 
$\SLC$-invariant algebraic   bidifferential operator  on $\TCP$.
\end{thm}

\begin{proof}   
We return to the proof of Proposition \ref{prop:4}.
Let $\bg_d=g_dD^d$ and $\bh_d=h_dD^d$, with $\bg_0=\bh_0=1$.
Writing out (\ref{eq:star=short}) termwise, we get, for $p\ge 1$, 
\begin{equation}\label{eq:Cp=long} 
 C_p(\phi,\psi)=\tsum_{i+j+k+m=p}\, \bg_mN_k(\bh_i\phi,\bh_j\psi)
\end{equation}
More succinctly,
$ C_p=\tsum_{i+j+k+m=p}\, \bg_mN_k(\bh_i\otimes\bh_j)$.

For $\la=\half$, the formula (\ref{eq:bh}) simplifies in that $[\half[d+1]]$ factors   
cancel out. Then  $h_d=U_dV_d\inv$ where 
$U_d=\invv[2^dd!]\prod_{i=1}^{[\half[d]]}(E'+i-\half)$ and
$V_d=\prod_{i=[\half[d+1]]}^{d-1}(E'+i)$.
Then   $\bh_d=U_dV_d\inv D^d=U_dD^dS_d\inv$. 
This is a   formal  relation,   valid for  $n$ odd  since then  $S_d$ is  
invertible. Similarly,  (\ref{eq:bg}) gives  
$\bg_d=T_d\inv\,F_dD^d$  where
$F_d=\invv[2^dd!]\prod_{i=[\half[d+1]]}^{d-1}(-E'-i-\half)$.
We put $U_0=F_0=1$.

We put $Z_p(\phi,\psi)=T_pC_p(S_p\phi,S_p\psi)$.
Let $T_{p;j}=T_p T_j\inv$ and  $S_{p;j}=S_j\inv S_p$.
Now (\ref{eq:Cp=long}) gives    $Z_p=\tsum_{i+j+k+m=p}Z^{mkij}$ where
\begin{equation}\label{eq:Zmkij=} 
\quad Z^{mkij}=T_{p;m}\, F_mD^m
N_k\left(U_iD^iS_{p;i}\otimes U_jD^jS_{p;j}\right)
\end{equation}
Each $Z^{mkij}$,
and so also their sum  $Z_p$,  is a bidifferential operator  on
$\TR$ with polynomial coefficients. I.e.,  
$Z_p$ lies in $\E\otimes_{\cP}\E$ where 
$\E=\C[u_i,\xi_j,\pdb{u_k},\pdb{\xi_l}]$  and $\cP=\C[u_i,\xi_j]$.

Now  $Z_p$ is invariant under $\slR$; this is clear since 
$T_p$, $C_p$ and $S_p$ are all  invariant.
It follows by projective geometry (as in \cite[\S8.1]{L-O}) that $Z_p$ extends 
uniquely to a global $\SLR$-invariant bidifferential operator on $\TRP$. 

We have $\{C_p(\phi,\psi)\,|\, \phi,\psi\in\RR\}$ $\to$
$\{Z_p(\phi,\psi)\,|\, \phi,\psi\in\RR\}$ $\to$
$\{Z_p(\phi,\psi)\,|\, \phi,\psi\in\A\}$ $\to$ 
$\{C_p(\phi,\psi)\,|\, \phi,\psi\in\A\}$  where the arrows 
indicate that one  set of values completely determines the next set.
The middle arrow follows because any bidifferential operator  on
$\TRP$ is completely determined by its values on $\RR$
(\cite[Lemma 5.1]{me:edq}).

Clearly $Z_p$ extends naturally (and uniquely) to an  algebraic
 differential operator  $\wt{Z}_p$   on $T^*\C^n$; 
this amounts to replacing our Darboux coordinates $u_i,\xi_j$
by their holomorphic counterparts $z_i,\ze_j$.
Then  $\wt{Z}_p$  is $\slC$-invariant and   (by projective geometry again)  
extends to  $\TCP$.
\end{proof}

Notice that this proof gives an explicit formula (in the coordinates $u_i,\xi_j$) for
$Z_p$.

\begin{remarks}\label{rems:5}   
(i) Suppose  $n$ is even.  Then this proof    still shows  
that the formula $Z_p(\phi,\psi)=T_pC_p(S_p\phi,S_p\psi)$ defines
an operator $Z_p$ in $\E\otimes_{\cP}\E$.
Then   (\ref{eq:Cp=}) is valid  as long as   $\phi$ and $\psi$ lie in 
$\A^*=\oplus_{d=[\half[p]]-\half[n]+1}^\infty\,\A^d$. We can
show that  all  the other results in Theorem \ref{thm:5} are still true, so that
(\ref{eq:Cp=}) determines $Z_p$ uniquely even for
$\phi,\psi\in\RR\cap\A^*$, 
$Z_p$ is an  $\SLR$-invariant bidifferential operator  on $\TRP$, etc.

(ii)   The maps  $\Q_{norm}$ and $\bh\sla$ 
are equivariant   with respect to only a  parabolic subgroup $P$ of $\SLR$,
even though their product $\Q\sla=\Q_{norm}\bh\sla$ is equivariant for $\SLR$.
Here $P$  is the subgroup of
the affine transformations of $\R^n$ (i.e., the one which fixes the
subspace $(u_0=0)$ in $\RPn$).    
Our formula (\ref{eq:Cp=})   is manifestly equivariant for  $\SLR$.
\end{remarks}

\section{Operators  $\Cla_p(\phi,\cdot)$ for  $\la=\half$}  
\label{sec_6}    

Next we  recover  part of  the results found for $\g=\fsl_{n+1}(\C)$
in \cite[Prop. 4.2.3]{A-B:exotic} and \cite[Thm. 6.3 and Cor. 8.2]{A-B:starmin}.
\begin{cor}\label{cor:6} 
Let $n\ge 1$.
For any momentum function $\mu^x$, $x\in\slC$, we have
\begin{equation}\label{eq:mars} 
C_2(\mu^x,\psi)=\frac{1}{E'(E'+1)}L^x(\psi), \qquad \psi\in\A
\end{equation}
where $L^x$ is an order $4$   differential operator  on 
$\TRP$.   

Neither $E'$ nor $E'+1$ left divides $L^x$ \textup{(}$x\neq 0$\textup{) }
over $T^*U$ for any open set $U$ in $\RPn$. Hence
$C_2(\mu^x,\cdot)$ is not a differential operator on  $T^*U$.  

Finally, $L^x$ extends uniquely to an  algebraic   differential operator 
on  $\TCP$. 
\end{cor}

\begin{proof} 
Suppose $n$ is odd. For   $\psi\in\A$, (\ref{eq:Cp=})  gives 
\begin{equation}\label{eq:venus} 
\textstyle  C_2(\mu^x,\psi)=
\frac{1}{E'+1}Z_2\left(\frac{1}{E'-1}\mu^x,\frac{1}{E'-1}\psi\right)=
\frac{2}{nE'(E'+1)}Z_2(\mu^x,\psi) 
\end{equation}   
The last equality follows because the operator $Z_2(\mu^x,\cdot)$ is graded of 
degree  $-1$.    

For $n$ even, (\ref{eq:venus}) is still true on account of  Remark \ref{rems:5}(i),
except in the case where
$n=2$ and $\psi\notin\oplus_{d=1}^\infty\,\A^d$. But if $\psi\in\A^0$ then
both  $C_2(\mu^x,\psi)$ and   $Z_2(\mu^x,\psi)$ vanish  for degree reasons and
so  the first and third expressions in  (\ref{eq:venus}) are still equal.
 
This proves  (\ref{eq:mars}), for all $n$, 
where $L^x=\frac{2}{n}Z_2(\mu^x,\cdot)$. 
Then $L^x$ extends to an  algebraic    differential operator  on $\TCP$; this follows
since both $Z_2$ and $\mu^x$ so extend. 

The  $L^x$, for $x\neq 0$, all have the same order. This follows because
the $L^x$, like the $\mu^x$,  transform in the adjoint \rep\ of $\SLC$.
We can choose  $\mu^x=\xi_m$ (the choice of
$m\in\{1,\dots,n\}$ is arbitrary). Let $L^{(m)}$ be the corresponding
operator  $L^x$. Using (\ref{eq:Cp=long}) we find after some calculation
\begin{equation}\label{eq:C2=} 
C_2(\cdot,\xi_m)=-\frac{1}{16}\frac{1}{E'(E'+1)}\xi_mD^2+
\frac{1}{8}\frac{1}{E'+1}\pdb{u_m}D   
\end{equation}
So  $L^{(m)}=-\frac{1}{16}(\xi_mD-2E'\pdb{u_m})D$. Clearly $L^{(m)}$ has 
order $4$. Using principal symbols, we see that $L^{(m)}$ has no left factors of 
the  form  $E'+c$ if $n\ge 2$. For $n=1$,  (\ref{eq:C2=})  gives
$L^{(m)}=\frac{1}{16}(E'+\half)\pdb[3]{u_1^2\del\xi_1}$, 
and so the only such  factor is   $E'+\half$.
\end{proof}

%We could prove the whole Corollary directly from (\ref{eq:C2=}), but it takes a
%little more work.

\begin{cor}\label{cor:last} 
Assume $n$ is odd and let $p\ge 1$.  If $\phi\in\A^d$ then
\begin{equation}\label{eq:last} 
C_p(\phi,\cdot)=
\invv[\prod_{i=1}^{[\frac{p}{2}]}\,(E'+i)(E'-i+p-d)]\, L^{\phi}_p 
\end{equation}
where $L^{\phi}_p$ is a differential operator on $\TRP$.
If $\phi\in\RR$, then  $L^{\phi}_p$ is an algebraic  
differential operator on $\TCP$.
\end{cor}
\begin{proof} 
This follows because $L^{\phi}_p=Z_p(S_p\inv\phi,\cdot)$.
\end{proof}  
 
%\begin{remark}\label{rem:SO}  
%We expect   similar  results for the   equivariant star
%products constructed in \cite{D-L-01}; it would be very nice to work this out. 
%There $\GR$ is the conformal group $SO(p+1,q+1)$ and $\XR$ is the manifold
%parameterizing  null lines through the origin in $\R^{p+1,q+1}$.  (The corresponding
%complex nilpotent orbit
%$\OO$ here is \emph{not} minimal, and so these cases are different from the $SO(n,\C)$
%cases treated in \cite{A-B}.)
%\end{remark}

\bibliographystyle{plain}

\end{document}